\documentclass[12pt,a4paper]{amsart}
\usepackage{amsfonts}
\usepackage{amssymb}
\usepackage{bbm}
\usepackage{comment}
\usepackage{enumitem}
\usepackage{dsfont}
\usepackage{amsthm}
\usepackage{amsmath}
\usepackage{commath}
\usepackage{url}
\usepackage{epsfig}
\usepackage{subfigure}
\usepackage{adjustbox}
\usepackage{mathtools}
\mathtoolsset{showonlyrefs,showmanualtags}

\makeatletter
\@namedef{subjclassname@2010}{%
  \textup{2010} Mathematics Subject Classification}
\makeatother

\usepackage[T1]{fontenc}

\newtheorem{theorem}{Theorem}
\newtheorem{corollary}[theorem]{Corollary}
\newtheorem{lem}[theorem]{Lemma}
\newtheorem{proposition}[theorem]{Proposition}

\newtheorem*{rmk}{Remark}

\theoremstyle{definition}

\newcommand{\al}{\alpha}
\newcommand{\ls}{\leq}
\newcommand{\gs}{\geq}

\def \bR {\mathbb R}

\title{Metric Density results for the value distribution of Sudler products}
\author{Manuel Hauke}
\address{Graz University of Technology, Institute of Analysis and Number Theory, Steyrergasse 30/II, 8010 Graz, Austria}
\email{hauke@math.tugraz.at}


\begin{document}

\begin{abstract}
    We study the value distribution of the Sudler product $P_N(\al) := \prod_{n=1}^{N}\lvert2\sin(\pi n \al)\rvert$ for Lebesgue-almost every irrational $\al$. We show that for every non-decreasing function \\\mbox{$\psi: (0,\infty) \to (0,\infty)$} with $\sum_{k=1}^{\infty} \frac{1}{\psi(k)} = \infty$, the 
    set $\{N \in \mathbb{N}: \log P_N(\al) \leq -\psi(\log N)\}$ has upper density $1$, which answers a question of Bence Borda.
    On the other hand, we prove that $\{N \in \mathbb{N}: \log P_N(\al) \geq \psi(\log N)\}$ has upper density at least $\frac{1}{2}$, with remarkable equality if $\liminf_{k \to \infty} \psi(k)/(k \log k) \geq C$ for some sufficiently large $C > 0$.
    \end{abstract}

\maketitle

\section{Introduction and statement of results}
For $\al \in \mathbb{R}$ and $N$ a natural number, the Sudler product is defined as
\begin{equation}\label{sudler}P_N(\al) := \prod_{r=1}^{N} 2 \left\lvert \sin (\pi r\al) \right\rvert.
\end{equation}
This product was first studied by Erd\"os and Szekeres \cite{erdos_szekeres}.
Later, Sudler products appeared in many different areas of mathematics that include, among others,
Zagier’s quantum modular forms and hyperbolic knots in algebraic topology \cite{quantum_invariants,bd2,zagier}, restricted partition functions \cite{sudler}, KAM theory \cite{kam} and Padé approximants \cite{lubinsky_pade}. Furthermore, they were used in the solution of the Ten Martini Problem \cite{aj}. Note that by $1$\,--\,periodicity of $P_N(\al)$ and the fact that $P_N(\al) = 0$ for rational $\al$ and $N$ sufficiently large,
it suffices to consider irrational numbers $\al \in [0,1]$.\\

In \cite{erdos_szekeres}, it was proven that
\begin{equation}\label{erd_sze}
\liminf_{N \to \infty} P_N(\al) = 0, \qquad \limsup_{N \to \infty} P_N(\al) = \infty
\end{equation}
holds for almost every $\al$, raising the question of whether this holds for all irrationals $\al$. 
Lubinsky \cite{lubinsky} showed that \eqref{erd_sze} remains true for all $\al$ that have unbounded partial quotients. On the other hand,
Grepstad, Kaltenböck and Neumüller showed in \cite{grepstad} that $\liminf_{N \to \infty} P_N(\phi) > 0$
for $\phi$ being the Golden Ratio, answering the question negatively. This counterexample was extended in \cite{tech_zaf,sudlerII} to certain quadratic irrationals that have only particularly small partial quotients.
For more results in this area, we refer the reader to \cite{grepstad_survey} and the references therein.\\

The asymptotic behaviour of the Sudler product depends delicately on the size of the partial quotients of $\al$. Since very much is known about the Diophantine properties for almost all irrationals, many results have been obtained in the metrical setting. 
Note that after taking logarithm, we see that $\log P_N(\al) = \sum_{r = 1}^N f(n \al)$ is a Birkhoff sum for the irrational rotation with 
$f(x) = \log\lvert 2 \sin(\pi x)\rvert$, having a logarithmic singularity. For a general
overview of Birkhoff sums in similar settings, we refer the reader to the survey \cite{birkhoff_surv}.
Lubinsky and Saff \cite{lub_saf}
proved that for almost all $\al$, we have
$\lim_{N \to \infty} \frac{\log P_N(\al)}{N} = 0$. Subsequently, Lubinsky \cite{lubinsky} improved this result and obtained
a divergence/convergence result as it is typical in metric Diophantine approximation: under a regularity condition
(see \cite{lubinsky} for the precise requirements), he showed that for a positive, non-decreasing function $\psi$ with
$\sum_{k =1}^{\infty} \frac{1}{\psi(k)} < \infty$, almost all $\al$ satisfy
\begin{equation}\label{conv_case}\lvert \log P_N(\al) \rvert \ll \psi(\log N)\end{equation}
(where $\ll$ denotes the usual Vinogradov symbol, see Section \ref{notation} for a proper definition).
On the other hand, if $\sum_{k =1}^{\infty} \frac{1}{\psi(k)} = \infty$,
then both inequalities
\begin{equation}\label{div_case}\log P_N(\al) \geq \psi(\log N), \qquad \log P_N(\al) \leq -\psi(\log N)\end{equation} hold for infinitely many $N$. These statements also follow from a more refined result obtained by Aistleitner and Borda \cite{quantum_invariants}, who showed that for all
$\al$ whose partial quotients fulfill \\\mbox{$(a_1 + \ldots + a_K)/K \to \infty$}, we have
\begin{equation}\label{max_sudler}\max_{0 \leq N < q_k} \log P_N(\al) = (V + o(1))(a_1 + \ldots + a_K) + \mathcal{O}\left(\frac{\log \max_{1 \leq \ell \leq K} a_{\ell}}{a_{K+1}}\right),\end{equation}
where $V = \int_{0}^{5/6} \log \lvert 2 \sin(\pi x) \rvert \,\mathrm{d}x \approx 0.1615$.
In a recent work, Borda \cite{borda} proved several results on the value distribution of Sudler products, both for 
badly approximable irrationals and for almost all $\al$. In the latter context, he improved \eqref{div_case}
in the sense that the inequalities in \eqref{div_case} both hold on a set of positive upper density.\\

\noindent {\bf Theorem A }(Borda, \cite[Theorem 6]{borda}). {\it Let $\psi$ be a non-decreasing, positive function such that $\sum_{k=1}^{\infty} \frac{1}{\psi(k)} = \infty$. Then for almost all $\al$, the sets
\begin{align}\label{first_set_borda}&\{N \in \mathbb{N}: \log P_N(\al) \geq \psi(\log N)\}\\
&\label{second_set_borda} \{N \in \mathbb{N}: \log P_N(\al) \leq -\psi(\log N)\}\end{align}
have upper density at least $\pi^2/(1440V^2) \approx 0.2627$, where 
$V = \int_{0}^{5/6} \log \lvert 2\sin(\pi x) \rvert \,\mathrm{d}x$.}\\

The proof relies on \eqref{max_sudler} and the variance estimate
\begin{equation}\label{variance_estim}\sqrt{\frac{1}{M}\sum_{N =1}^M \log^2 P_N(\al)} = \left(\frac{\pi}{\sqrt{720}V} + o(1)\right)\max_{0 \leq N < M} \log P_N(\al), \end{equation}
which is shown to hold for infinitely many $M \in \mathbb{N}$. 
Additionally, Borda makes use of
the ``reflection principle'' of Sudler products, which will also play a main role in this paper. This principle was observed by \cite{tech_zaf} and used in the subsequent literature on Sudler products several times. We state it here in the form of \cite[Propositions 2 and 3]{quantum_invariants}:
for any irrational $\al$ and $0 \leq N < q_K$ (where $q_K$ denotes the denominator of the $k$\,--\,th convergent of $\alpha$, see Section \ref{cont_frac} for a proper definition), we have
\begin{equation}\label{reflection_princ}\log P_N(\al) + \log P_{q_K-N-1}(\al) = \log q_K + \mathcal{O}\left(\frac{1 + \log \max_{\ell \leq K} a_{\ell}}{a_{K+1}}\right).\end{equation}
In particular, \eqref{reflection_princ} implies that for almost all $\al$, the values $\log P_N(\al),\; N = 1,\ldots, q_K,$ distribute symmetrically around the center $\log q_K$, which is however of negligible order for almost all $\al$. Hence, the numbers $1 \leq N < q_K$
lie approximately as often in \eqref{first_set_borda} as in \eqref{second_set_borda}.
Borda remarked in \cite{borda} that the estimate on the upper density in Theorem A is probably not optimal, saying that it might be possible that the union of \eqref{first_set_borda} and \eqref{second_set_borda} has upper density $1$.
Here we prove something even stronger: we show that already
\eqref{second_set_borda} on its own has upper density $1$.

\begin{theorem}\label{thm1}
Let $\psi$ be a non-decreasing, positive function such that $\sum_{k=1}^{\infty} \frac{1}{\psi(k)} = \infty$. Then for almost every $\al$, the set
\[\{N \in \mathbb{N}: \log P_N(\al) \leq -\psi(\log N)\}\]
has upper density $1$.
\end{theorem}

The symmetry around the negligible center $\log q_k$ discussed above leads to the belief that \eqref{first_set_borda} has the same upper density 
than \eqref{second_set_borda}.
Surprisingly, this turns out to be wrong: we prove that if $\psi$ is as in Theorem \ref{thm1} and additionally fulfills a certain regularity condition,
\eqref{first_set_borda} has upper density $1/2$ for almost every $\al$.

\begin{theorem}\label{thm2}
Let $\psi$ be a non-decreasing, positive function such that $\sum_{k=1}^{\infty} \frac{1}{\psi(k)} = \infty$.
Then for almost every $\al$, the set

\[\{N \in \mathbb{N}: \log P_N(\al) \geq \psi(\log N)\}\]
has upper density at least $1/2$, with equality if $\liminf\limits_{k \to \infty} \frac{\psi(k)}{k \log k} \geq C$ for some absolute constant $C > 0$.
\end{theorem}

\subsection*{Remarks on Theorems \ref{thm1} and \ref{thm2} and further research}
\begin{itemize}
\item Note that the divergence criterion of $\sum_{k=1}^{\infty}\frac{1}{\psi(k)}$ is invariant under multiplication with constant factors. Therefore, it suffices to show Theorems \ref{thm1} and the first part of Theorem \ref{thm2} for the sets \eqref{first_set_borda} and \eqref{second_set_borda} with $\psi (\log N)$ substituted with $C_1\cdot \psi(C_2 \log N)$, where $C_1,C_2 > 0$ are arbitrary constants. We will make use of this fact several times in the subsequent proofs without explicitly stating it.\\

\item By \eqref{conv_case}, we see that the assumption $\sum_{k=1}^{\infty} \frac{1}{\psi(k)} = \infty$ is essential, as otherwise the
upper density is trivially zero.
 Note that also ``upper density'' cannot be replaced by ``lower density'': for $\psi(k) \geq \left(12V/\pi^2 + \varepsilon\right)k \log k$, where $V$ is the constant from Theorem A, even the union of \eqref{first_set_borda} and \eqref{second_set_borda} has lower density zero (see \cite[Theorem 7]{borda}). It is interesting to find the minimal growth rate of $\psi$ such that the sets \eqref{first_set_borda}, \eqref{second_set_borda} or their union have non-zero lower density.
 \\

\item Note that even in the case when the regularity condition $\liminf\limits_{k \to \infty} \frac{\psi(k)}{k \log k} \geq C$ is not satisfied, Theorem 2 gives an improved lower bound in comparison to Theorem A.
Our approach relies on the fact that for almost every irrational, the trimmed sum of its first $k$ partial quotients is bounded from above by $k \log k$, with the largest partial quotient dominating the sum infinitely often. Therefore, we only need to control the Ostrowski coefficient of the largest partial quotient (see Section \ref{heuristic} for an overview).
It remains open how far the regularity condition from Theorem \ref{thm2} can be relaxed such that the upper density of \eqref{first_set_borda} is still $1/2$ for almost every $\al$. 
Below we show that $\psi$ has to fulfill $\psi(k) \geq (1/2 - \varepsilon)k$ infinitely often for arbitrary small $\varepsilon> 0$. This can be deduced in the following way from \cite[Theorem 9]{borda}:
the theorem states (among other results) that for any $t \geq 0$,

\begin{align}\lim_{M \to \infty}
\lambda&\Bigg(\Bigg\{ \alpha \in [0,1] : \frac{10\pi}{M \log^2 M}
\sum_{N =1}^M \left(\log P_N(\al) - \frac{1}{2}\log M\right)^2 \leq t\Bigg\}\Bigg)
\\= &\int_{0}^t \frac{e^{-1/(2x)}}{\sqrt{2\pi}x^{3/2}} \,\mathrm{d}x =: c(t),\end{align}
where $\lambda$ denotes the $1$\,--\,dimensional Lebesgue measure.
By Chebyshev's inequality, we obtain that for any $\varepsilon, y > 0$,

\begin{align}
\liminf_{M \to \infty} 
\lambda&\Bigg(\Bigg\{  \alpha \in [0,1] : \frac{1}{M} \#\left\{ 1 \leq N \leq M:
\left\lvert \log P_N(\al) - \frac{1}{2}\log M \right\rvert \geq \varepsilon\log M \right\} \leq y \Bigg\}\Bigg) 
\\\geq& \;c\;\!(10\pi \varepsilon^2 y).
\end{align}
Applying Fatou's Lemma, we get that on a set of measure of at least $ c\;\!(10\pi \varepsilon^2 y) > 0$,
\[
\frac{1}{M} \#\left\{1 \leq N \leq M:
\left\lvert \log P_N(\al) - \frac{1}{2}\log M \right\rvert \geq \varepsilon\log M\right\} \leq y
\]
holds for infinitely many $M$. This implies that the upper density of 
\[
\left\{N \in \mathbb{N}: \log P_N(\al) > \left(\frac{1}{2}-\varepsilon\right)\log N \right\}
\]
is bounded from below by $1-y$, so choosing $y < \frac{1}{2}$, we can deduce that for $\psi(k) \leq (1/2 - \varepsilon)k$, the upper density of \eqref{first_set_borda} being $1/2$ fails to hold on a set of positive measure.
However, it remains open whether having $\psi(k) \geq \frac{k}{2}$ is already sufficient to deduce upper density $1/2$ for almost all $\al$.\par
Similarly, it is interesting if there is some threshold function where the upper density of the set in \eqref{first_set_borda} jumps from $1/2$ to $1$ for almost every $\al$ (and if so, how fast does this function grow?), or if the value of the upper density attains a fixed number strictly between $1/2$ and $1$ for certain functions $\psi$ and almost every irrational.
\end{itemize}

\section{Notation and preliminary results}

\subsection{Notation}\label{notation} Given two functions $f,g:(0,\infty)\to \mathbb{R},$ we write $f(x) = \mathcal{O}\left(g(x)\right)$ or $f \ll g$, when  \vspace{-2mm}
$\limsup_{x\to\infty} \frac{f(x)}{g(x)} < \infty$ and $f(x) = o\left(g(x)\right)$,
when $\limsup_{x\to\infty} \frac{f(x)}{g(x)} = 0$.
If $f \ll g$ and $g \ll f$, we write $f \asymp g$ and $f \sim g$ for $\lim_{x\to\infty}\frac{f(x)}{g(x)} = 1$.
Given a real number $x\in \bR,$ we write $\|x\|=\min\{|x-k|: k\in\mathbb{Z}\}$ for the distance of $x$ from its nearest integer.

\subsection{Continued fractions}\label{cont_frac}
In this subsection, we shortly recall all necessary facts about the theory on continued fraction that are used to prove Theorems \ref{thm1} and \ref{thm2}. For a more detailed introduction, we refer the reader to the classical literature, e.g. \cite{all_shall, rock_sz, schmidt}. Every irrational $\al$ has a unique infinite continued fraction expansion $[a_0;a_1,...]$ with convergents $p_k/q_k = [a_0;a_1,...,a_k]$ that fulfill the recursions 
\begin{equation}\label{rec_converg}p_{k+1} = p_{k+1}(\al) = a_{k+1}(\al)p_k + p_{k-1}, \quad q_{k+1} = q_{k+1}(\al) = a_{k+1}(\al)q_k + q_{k-1}.\end{equation}
For shorter notation, we will just write $p_k,q_k,a_k$, although these entities depend on $\al$.
We know that $p_k/q_k$ approximates $\al$ very well, which leads to the following well-known inequalities for $k \gs 1$:

\begin{equation}\label{approx_cf}\frac{1}{q_{k+1} + q_k} \leq \delta_k := \lVert q_k \al \rVert = \lvert q_k\al - p_k \rvert \ls \frac{1}{q_{k+1}},\end{equation}
from where we can deduce that
\begin{equation}\label{qkqkalpha}
    \frac{1}{a_{k+1}+2}\ls q_k\delta_k \leq \frac{q_k}{q_{k+1}} \ls \frac{1}{a_{k+1}}.
\end{equation}
Using \eqref{rec_converg}, we obtain that
\begin{equation}\label{q_kal_in_other}
    a_{k+1}\delta_k = \delta_{k-1} - \delta_{k+1}.
\end{equation}
Fixing an irrational $\al = [a_0;a_1,...]$, the Ostrowski expansion of a non-negative integer $N$ is the unique representation

\begin{equation}\label{ostrowski}N = \sum_{\ell = 0}^{K-1} b_{\ell}q_{\ell} \quad \text{ where }
\quad b_{K} \neq 0,  \quad 0 \ls b_{\ell} \ls a_{\ell +1}, \quad b_0 < a_1,
\end{equation}
with the additional rule that
$b_{\ell-1} = 0$ whenever $b_{\ell} = a_{\ell+1}$.\vspace{5mm}\par

\subsection*{Metrical results}
Much is known about the almost sure behavior of continued fraction coefficients and convergents.
Below we state all known properties of almost every $\al$ that are used during the proofs of Theorems \ref{thm1} and \ref{thm2}. 

\begin{itemize}
    \item (Bernstein \cite{bernstein}): For any monotonically non-decreasing function $\psi: [1,\infty) \to [1,\infty)$, we have
    \begin{equation}
        \label{bernstein}
        \#\left\{k \in \mathbb{N}: a_k > \psi(k)\right\} \text{ is } \begin{cases} \text{ infinite \hspace{5mm}if } \sum_{k = 0}^{\infty} \frac{1}{\psi(k)} = \infty
        \\ \text{ finite \hspace{7mm} if }\sum_{k = 0}^{\infty} \frac{1}{\psi(k)} < \infty.\end{cases}
    \end{equation}
\item (Diamond and Vaaler \cite{diamond_vaaler}):\begin{equation}\label{trimmed_sum_eq}\sum_{\ell \leq K} a_{\ell} - \max\limits_{\ell \leq K} a_{\ell} \sim \frac{K \log K}{\log 2}, \quad K \to \infty.\end{equation} 
\item (Khintchine and L\'{e}vy, see e.g. \cite[Chapter 5, §9, Theorem 1]{rock_sz}):
\begin{equation}\label{size_of_q_k}\log q_k \sim \tfrac{\pi^2}{12 \log 2} k \text{ as } k \to \infty.\end{equation}
\end{itemize}

Combining \eqref{bernstein} and \eqref{trimmed_sum_eq}, the following corollary follows immediately.

\begin{corollary}\label{needed_metric_results}
Let $\psi$ be a non-decreasing, positive function such that $\sum_{k=1}^{\infty} \frac{1}{\psi(k)} = \infty$.
Then for almost every $\al$, there exist infinitely many $K \in \mathbb{N}$ such that the following hold.
\begin{itemize}
    \item[a)] $\psi(K) < a_K < K^2 $.
    \item[b)] $\sum_{\ell =1}^{K-1} a_{\ell} \ll K \log K$
    with an absolute implied constant.%, independent of $K$.
\end{itemize}
\end{corollary}

\section{Heuristic behind the proofs}\label{heuristic}

We start by sketching the heuristic idea behind the proof of Theorems \ref{thm1} and \ref{thm2}.
This can be compared with \cite[Section 2.1]{aist_borda}. Starting with Theorem \ref{thm1}, note that we can assume without loss of generality that $\psi(k)/(k \log k) \to \infty$, since this implies the statement also for slower-growing $\psi$.
Let $\psi$ and $K$ be as in Corollary \ref{needed_metric_results} and let $N < q_K$ be arbitrary with Ostrowski expansion
$N = \sum_{\ell=0}^{K-1} b_{\ell}q_{\ell}$.
We use the usual decomposition of $P_N(\al)$ into certain shifted Sudler products. This approach was first used in the special case for $\al$ being the Golden Ratio in \cite{grepstad} and made more explicit and general in subsequent works in this area, e.g. \cite{quantum_invariants,tech_zaf,grepstad_survey,sudlerII,myother}. Defining
\[P_N(\al,x) := \prod_{n = 1}^N \lvert 2 \sin(\pi(n\al +x ))\rvert, \quad \al,x \in \mathbb{R},\]
and
\begin{equation}\label{def_of_eps}
\varepsilon_{\ell}(N) := q_{\ell}\sum_{k = \ell+1}^{K-1}(-1)^{k+\ell}b_{k}\delta_k,
\end{equation}
we can deduce (see \cite[Lemma 2]{quantum_invariants}) that
\begin{equation}\label{shifted_decomp}P_N(\al) = \prod_{\ell = 0}^{K-1}\prod_{b = 0}^{b_{\ell}-1}P_{q_{\ell}}\left(\al,(-1)^{\ell}(bq_{\ell}\delta_{\ell} + \varepsilon_{\ell}(N))/q_{\ell}\right).
\end{equation}
Ignoring first the contribution of the numbers $\varepsilon_{\ell}(N)$, and using the approximation\\
$P_{q_{\ell}}\left(\al,(-1)^{\ell}x/q_{\ell}\right) \approx \lvert 2 \sin(\pi x)\rvert$ elaborated later, we see that

\begin{align}
    \log P_N(\al) &\approx \sum_{b = 1}^{b_{K-1} -1}\log \left\lvert 2\sin\left(\pi bq_{K-1}\delta_{K-1}\right) \right\rvert
    + \sum_{\ell=0}^{K-2}\sum_{b = 1}^{b_{\ell} -1}\log \left\lvert 2\sin\left(\pi bq_{\ell}\delta_{\ell}\right) \right\rvert\\&\approx a_K\int_{0}^{b_{K-1}/a_{K}} \log \left\lvert 2\sin(\pi x)\right\rvert \,\mathrm{d}x + \sum_{\ell=0}^{K-2}
    a_{\ell+1}\int_{0}^{b_{\ell}/a_{\ell+1}} \log \left\lvert 2\sin(\pi x)\right\rvert \,\mathrm{d}x.
\end{align}
By the choice of $K$ as in Corollary \ref{needed_metric_results}, the value $a_K$ dominates the sum $\sum_{\ell = 0}^{K-1} a_{\ell}$. So using $\log \lvert 2\sin(\pi x)\rvert \leq \log(2)$ and assuming that 
\begin{equation}\label{cruc_integral}\int_{0}^{b_{K-1}/a_{K}} \log \lvert 2\sin(\pi x)\rvert \,\mathrm{d}x\end{equation} is bounded away from $0$, we have that
$\log P_N(\al) \ll -a_K$, provided that the integral in \eqref{cruc_integral} is negative. It is easy to see that this is the case if and only if $b_{K-1}/a_{K} < \frac{1}{2}$, which leads to
\begin{equation}\label{less_heur}\log P_N(\al) \ll -\psi(K)\end{equation} for $b_{K-1}/a_{K} < 1/2 - \varepsilon$.
As almost all numbers $N < \left\lfloor\frac{q_k}{2}\right\rfloor$ fulfill 

\[\log N \asymp \log q_K \underset{\eqref{size_of_q_k}}\asymp K,\] 
\eqref{less_heur} is equivalent to $\log P_N(\al) \ll -\psi(\log N)$ for most $N$, which implies Theorem \ref{thm1}.\\

By the same reasoning, we can immediately deduce that at least $50\%$ 
of all numbers $N < q_K$ fulfill \eqref{less_heur}.
Using the reflection principle, we see that also
\begin{equation}\label{more_heur}\log P_N(\al) \gg \psi(K)\end{equation} is fulfilled for 
about $50\%$ of all numbers $N < q_K$,
hence the first part of Theorem \ref{thm2} follows immediately.
For the equality in case $\liminf_{k \to \infty} \psi(k)/(k \log k) \geq C$, we fix some integer $q_{K-1} \leq M < q_K$ 
(this $K$ does not fulfill in general the properties of Corollary \ref{needed_metric_results}),
and show that
asymptotically, at most 50\% of all $N < M$ can fulfill $\log P_N(\al) \gg \psi(K)$.
Defining \mbox{$a_{{\ell_0}} = \max_{\ell \leq K} a_{\ell}$}, we can argue similar to before that for $C$ sufficiently large and $\log N \gg \log q_K$,
\begin{align*}\log P_N(\al) &\lesssim a_{{\ell_0}}\int_{0}^{b_{{\ell_0}-1}/a_{{\ell_0}}} \log \lvert 2\sin(\pi x)\rvert \,\mathrm{d}x + \mathcal{O}\left(\sum_{k\neq {\ell_0}}^{K}a_k\right)\\
&\leq a_{{\ell_0}}\int_{0}^{b_{{\ell_0}-1}/a_{{\ell_0}}} \log \lvert 2\sin(\pi x)\rvert \,\mathrm{d}x + \frac{\psi(\log N)}{2}.\end{align*}
So in order to fulfill $\log P_N(\al) \geq \psi(\log N)$, we have the necessary condition
\begin{equation}\label{cond_pos_int}\int_{0}^{b_{{\ell_0}-1}/a_{{\ell_0}}} \log \lvert 2\sin(\pi x)\rvert > 0,\end{equation} or equivalently,
$b_{{\ell_0}-1}(N)/a_{{\ell_0}} > 1/2$, which can be seen to be fulfilled by at most $50$\% of all $N < M$. 
Hence, no matter how we choose $M \in \mathbb{N}$, at most half the numbers $N < M$ fulfill \eqref{cond_pos_int}, so the upper density of \eqref{first_set_borda} cannot exceed $1/2$.\\

The punchline why the upper densities of \eqref{first_set_borda} and \eqref{second_set_borda} differ is the following: on the full period $1 \leq N \leq q_K$, there are about as many elements in \eqref{first_set_borda} as in \eqref{second_set_borda}, and for $a_K$ being large, almost all elements are in one of those sets.
The criterion whether $N$ is in \eqref{first_set_borda} or in \eqref{second_set_borda} is (almost) equivalent $b_{K-1}(N) > a_{K}/2$ or not.
As $b_{K-1}$ is the most significant coefficient for the size of $N$ (since $b_{K-1}(M) < b_{K-1}(N)$ implies $M < N$), we see that
all elements in \eqref{second_set_borda} appear before the elements in \eqref{first_set_borda}, causing the asymmetric result.

\begin{rmk}
Note that all estimates in this paper only consider upper bounds. This makes
the analysis much easier since we can ignore the singularities of the function $\log\lvert 2\sin(\pi x)\rvert$ at $x = 0$ or $x=1$, as we trivially bound $\log\lvert 2\sin(\pi x)\rvert \leq \log(2)$ from above. The reflection principle provides the tool to use the upper bounds also to achieve Theorem \ref{thm2}, without having to consider that singularities.
\end{rmk}

\section{Proof of the theorems}
\subsection{Preparatory results for the approximation errors}
In this section, we discuss the actual errors that are made by comparing $\log P_N(\al)$ with
$a_K\int_{0}^{b_{K-1}/a_{K}} \log \lvert 2\sin(\pi x)\rvert \,\mathrm{d}x$ (see Lemma \ref{upper_bound_key_lem}). The first step in this direction is done by 
\cite[Proposition 12]{aist_borda}. For the convenience of the reader, we state it below as Proposition \ref{prop12}.

\begin{proposition}\label{prop12}
Let $N = \sum_{\ell = 0}^{K-1}b_{\ell}q_{\ell}$ be the Ostrowski expansion of a non-negative integer and $\varepsilon_{\ell}(N)$ as in \eqref{def_of_eps}. There exists a universal constant $C > 0$ such that for any $\ell \geq 1$ with $b_{\ell} \geq 1$,
we have

\begin{align*}
\sum_{b = 0}^{b_{\ell}-1}\log P_{q_{\ell}}(\al,(-1)^{\ell}(bq_{\ell}\delta_{\ell} + \varepsilon_{\ell}(N))/q_{\ell}) \leq &\sum_{b = 1}^{b_{\ell}-1} \log \lvert 2 \sin(\pi(bq_{\ell}\delta_{\ell} + \varepsilon_{\ell}(N)))\rvert \\+ &\sum_{b = 0}^{b_{\ell}-1} V_{\ell}(bq_{\ell}\delta_{\ell} + \varepsilon_{\ell}(N))
\nonumber \\+ &\nonumber\log(2\pi(b_{\ell}q_{\ell}\delta_{\ell} + \varepsilon_{\ell}(N)) + \frac{C}{a_{{\ell}+1}q_{\ell}},
\end{align*}
where 
\begin{equation}\label{def_V_k}V_{\ell}(x) := \sum_{n = 1}^{q_{\ell}-1} \sin(\pi n \delta_{\ell}/q_{\ell}) \cot\left(\pi \frac{n (-1)^{\ell}p_{\ell} + x}{q_{\ell}}\right)\end{equation} denotes a modified cotangent sum.
\end{proposition}

We see that we need to find upper bounds on the modified cotangent sums $V_{\ell}$. This is done by the following variant of \cite[Lemma 8]{aist_borda}.

\begin{lem}\label{V_k_lem}
Let $1 \leq k \leq K-1$, $a_{{\ell_0}} = \max_{1 \leq \ell \leq K} a_{\ell}$, $x \in (-1,1)$ and $V_k$ as in \eqref{def_V_k}. Then the following statements hold.
\begin{itemize}
    \item[(i)] \[V_k'(x) < 0, \qquad \lvert V_k(0)\rvert \ll \frac{1+ \log a_{{\ell_0}}}{a_{k+1}}.\]
    \item[(ii)] \[\lvert V_k(x)\rvert \ll \log a_{{\ell_0}} + \frac{1}{1-\lvert x \rvert},\]
    with the implied constants independent of $x$ and $k$.
\end{itemize}
\end{lem}

\begin{proof}
The statements in (i) are proven in \cite[Lemma 8]{aist_borda}. For (ii), we use the estimate 
$\lvert V_k'(x)\rvert \ll \frac{1}{(1 - \lvert x \rvert)^2}$, which is also shown in \cite{aist_borda}. The result now follows immediately after integration.
\end{proof}

Next, we turn our attention to controlling the size of the perturbations $\varepsilon_{\ell}(N)$.
It is easy to see that $-1 < \varepsilon_{\ell}(N) < 1$ for any $1 \leq \ell \leq K-1$. By Lemma \ref{V_k_lem}, we see that 
the error made by  $V_{\ell}(bq_{\ell}\delta_{\ell} + \varepsilon_{\ell}(N))$ is particularly large
when its argument is close to its singularities at $-1$ and $1$. The following proposition aims to bound the arguments away from those singularities and to show that the perturbation $\varepsilon_{\ell}(N)$ is small if $a_{\ell +1}$ is large, which will be the case in the main term (see Section \ref{heuristic}).
\begin{proposition}\label{size_of_eps}
Let $\varepsilon_{\ell}(N)$ be defined as in \eqref{def_of_eps} and $b_{\ell} \geq 1$. Then we have the following inequalities:
\begin{itemize}
    \item[(i)]
\begin{equation}\label{in_part_eps}
    - \frac{1}{a_{{\ell}+1}} \leq -q_{\ell}\delta_{\ell} \leq \varepsilon_{\ell}(N) \leq \frac{1}{a_{{\ell}+1}}.
\end{equation}
\item[(ii)]
\begin{equation}\label{bound_from-1_1}1 - \rvert\varepsilon_{\ell}(N)\rvert \gg \frac{1}{a_{\ell +2}}.\end{equation}
If $b_{{\ell}+1} \leq \frac{a_{{\ell}+2}}{2}$, then
\begin{equation}\label{bound_from-1_2}1 - \rvert\varepsilon_{\ell}(N)\rvert \gg 1,\end{equation}
with the implied constants being absolute.
\end{itemize}
\end{proposition}

\begin{proof}
We argue similarly to \cite[Lemma 3]{quantum_invariants}. By definition of $\varepsilon_{\ell}(N)$ and \eqref{q_kal_in_other}, we obtain
\begin{align*}
    \varepsilon_{\ell}(N) = q_{\ell}\sum_{k = \ell+1}^{K-1}(-1)^{k+\ell}b_{k}\delta_k
    &\leq q_{\ell}(a_{\ell+3}\delta_{\ell +2} + a_{\ell+5}\delta_{\ell +4} + \ldots)
    \\&= q_{\ell}\left((\delta_{\ell +1} - \delta_{\ell +3}) + (\delta_{\ell +3} - \delta_{\ell +5}) + \ldots \right)
    \\&= q_{\ell}\delta_{\ell +1}
    \leq \frac{q_{\ell}}{q_{\ell +2}} \leq \frac{1}{2},
\end{align*}
where we used \eqref{approx_cf} in the last line.
Similarly, we get
\begin{align*}
    \varepsilon_{\ell}(N) &\geq -q_{\ell}(b_{\ell+1}\delta_{\ell +1} + a_{\ell+4}\delta_{\ell +3} + \ldots)
    \\&= -q_{\ell}\left((b_{\ell+1} - a_{\ell +2})\delta_{\ell +1} + (\delta_{\ell} - \delta_{\ell +2}) + (\delta_{\ell +2} - \delta_{\ell +4}) + \ldots \right)
    \\&= -q_{\ell}\left(\delta_{\ell} - (b_{\ell+1} - a_{\ell +2})\delta_{\ell +1}\right).
\end{align*}
As $b_{\ell} \geq 1$ implies $b_{\ell +1} \leq a_{\ell +2}-1$, combining these bounds leads to
    \begin{equation}\label{ineq_chain_eps} -1 < -q_{\ell}\delta_{\ell} + q_{\ell} \delta_{{\ell} +1} \leq
    -q_{\ell}\delta_{\ell} + q_{\ell}(a_{{\ell}+2}-b_{{\ell}+1}) \delta_{{\ell} +1} \leq \varepsilon_{\ell}(N) \leq q_{\ell} \delta_{{\ell} +1}
    \leq \frac{1}{2}.\end{equation}
    \begin{itemize}
    \item[(i):]
    As $\delta_{{\ell} +1} \leq \delta_{\ell}$, \eqref{in_part_eps} follows immediately from \eqref{qkqkalpha} and \eqref{ineq_chain_eps}.\\
    \item[(ii):]
    By \eqref{ineq_chain_eps}, we have $\varepsilon_{\ell}(N) < \frac{1}{2}$, so it suffices to find lower bounds for
    $\varepsilon_{\ell}(N)$.
    Using \eqref{approx_cf} and $q_{{\ell}+1} \leq 2a_{{\ell}+1}q_{{\ell}}$, we get
    
    \[q_{\ell}\delta_{{\ell} +1}
    \geq \frac{q_{\ell}}{q_{{\ell}+2} + q_{{\ell}+1}} \geq \frac{q_{\ell}}{3a_{{\ell}+2}q_{{\ell}+1}} \geq \frac{1}{6a_{{\ell}+2}a_{{\ell}+1}}.\]
    Applying \eqref{qkqkalpha}, we get
    \[-q_{\ell}\delta_{\ell} + (a_{{\ell}+2}-b_{{\ell}+1})q_{\ell}\delta_{{\ell}+1} \geq \frac{1}{a_{{\ell}+1}}\left(-1 + \frac{a_{{\ell}+2}-b_{{\ell}+1}}{6a_{{\ell}+2}}\right),\]
    which in view of \eqref{ineq_chain_eps} finishes the proof.
    \end{itemize}
    \end{proof}

The following lemma combines the preparatory results from above. It contains the main ingredients to the proof of both Theorems \ref{thm1} and \ref{thm2}.
\begin{lem}\label{upper_bound_key_lem}
Let  $N < q_K$ with Ostrowski expansion $\sum_{\ell=0}^{K-1} b_{\ell}q_{\ell}$ and
let $1 \leq {\ell_0} \leq K$ be such that $a_{{\ell_0}} = \max_{\ell \leq K} a_{\ell} \geq 2$.
Assume that $b_{{\ell_0}-1} \leq \frac{a_{{\ell_0}}}{2} \leq \frac{K^2}{2}$
and \begin{equation}\label{last_dominates_sum}\sum_{\substack{k =1,\\ k \neq {\ell_0}}}^{K} a_k \ll K \log K.\end{equation}
Then we have

\begin{equation}\label{upper_bound_eq}
    \log P_N(\al) \leq \sum_{b = 1}^{b_{{\ell_0}-1} -1}\log \left\lvert 2\sin\left(\pi bq_{{\ell_0}-1}\delta_{\ell_0-1}
    + \varepsilon_{{\ell_0}-1}(N)\right) \right\rvert
    + \mathcal{O}\left(K \log K\right).
\end{equation}
\end{lem}

\begin{proof}[Proof of Lemma \ref{upper_bound_key_lem}]
Using the decomposition into shifted Sudler products from \eqref{shifted_decomp}, we have
\begin{equation}
    \log P_N(\al) = \sum_{k=0}^{K-1}\sum_{b = 0}^{b_k -1}\log P_{q_k}\left(\al, (-1)^k(bq_k\delta_k 
    + \varepsilon_k(N))/q_k\right).
\end{equation}
Next, we apply Proposition \ref{prop12}
for every $1 \leq k \leq K-1$ with $b_k \neq 0$ and obtain for some $C > 0$ that

\begin{align*}
\log P_N(\al) \leq 
\sum_{k = 1}^{K-1}\Bigg(&\sum_{b = 1}^{b_k-1} \log \left\lvert 2 \sin(\pi(bq_k\delta_k + \varepsilon_k(N)))\right\rvert \\+ &\sum_{b = 0}^{b_k-1} V_{k}(bq_k\delta_k + \varepsilon_k(N))
\nonumber \\+ &\nonumber\log(2\pi(b_kq_k\delta_k + \varepsilon_k(N))) + \frac{C}{a_{k+1}q_k}\Bigg).
\end{align*}
Applying rough bounds on the arguments of the logarithms and using \eqref{last_dominates_sum} leads to
\begin{align}\log P_N(\al) 
\leq \hspace{2mm}&\sum_{b = 1}^{b_{{\ell_0}-1} -1}\log \left\lvert2\sin\left(\pi bq_{{\ell_0}-1}\delta_{\ell_0-1}+ \varepsilon_{{\ell_0}-1}(N)
    \right)\right\rvert \\+\hspace{2mm}& \hspace{2mm}\sum_{k = 1}^{K-1}\sum_{b = 0}^{b_k-1} V_{k}\left(bq_k\delta_k + \varepsilon_k(N)\right) + \mathcal{O}\left(K \log K\right).
\end{align}
By Proposition \ref{size_of_eps}\,(i), we see that $b \geq 1$ implies that
$bq_k\delta_k + \varepsilon_k(N) \geq 0$. So Lemma \ref{V_k_lem}\,(i) and $a_{{\ell}_0} \leq K^2$ lead to
\[\sum_{k = 1}^{K-1}\sum_{b = 1}^{b_k-1} V_{k}(bq_k\delta_k + \varepsilon_k(N)) \ll
\sum_{k = 1}^{K-1}\frac{b_k}{a_{k+1}}\log a_{{\ell_0}}
\ll K \log K.\]
For $1 \leq k\neq {\ell_0}-2 \leq K-2$, we use \eqref{bound_from-1_1} to obtain 
\begin{equation}\label{cotangent_0_1}\frac{1}{1 - \lvert \varepsilon_k(N) \rvert} \ll a_{k+2}.\end{equation}
For $k = {\ell_0}-2$, we observe that $b_{{\ell_0}-1} \leq \frac{a_{{\ell_0}}}{2}$, hence we have by \eqref{bound_from-1_2} that
\begin{equation}\label{cotangent_0_2}\frac{1}{1 - \lvert \varepsilon_{{\ell_0}-2}(N) \rvert} \ll 1.\end{equation}
For $k = {\ell_0}-1$, we apply Proposition \ref{size_of_eps}\,(i) to obtain $\lvert\varepsilon_{{\ell_0}-1}(N)\rvert \leq \frac{1}{2}$,
and from the definition of $\varepsilon_{\ell}(N)$, we can follow that $\varepsilon_{K-1}(N) = 0$. Combining these observations with \eqref{cotangent_0_1} and \eqref{cotangent_0_2} yields
\[\sum_{k = 1}^{K-1}V_{k}(\varepsilon_k(N)) \ll K \log K,\]
where we used  \eqref{last_dominates_sum} once more.
This finishes the proof.
\end{proof}

\subsection{Proof of Theorem \ref{thm1}}
We can assume without loss of generality that $\lim_{k \to \infty} \psi(k)/(k \log k) = \infty$, as showing this will imply the statement of Theorem \ref{thm1} also for slower growing $\psi$.
Applying Corollary \ref{needed_metric_results}, we know that there exist infinitely many $K$ such that
\begin{equation}\label{cor_3_prop}\psi(K) < a_K < K^2, \quad \sum_{k =1}^{K-1} a_k \ll K \log K.\end{equation}
Fixing an arbitrary small $\delta >0$, we define for every $K \geq 1$ that fulfills \eqref{cor_3_prop}, the set
\[M_K = M_K(\delta) := \left\{1 \leq N \leq \left\lfloor \frac{q_K}{2}\right\rfloor: \delta a_K \leq b_{K-1}(N) \leq \left(\frac{1}{2}-\delta\right)a_K\right\}.\]
Choosing $K$ sufficiently large, we have by \eqref{size_of_q_k}  that for all $N \in M_K$,

\begin{equation}\label{same_to_const}\psi(\log N)  \asymp \psi(\log M_K) \asymp \psi(K).\end{equation}
As $\# M_K(\delta) / \left\lfloor \frac{q_K}{2}\right\rfloor \underset{\delta \to 0}{\to} 1$, it suffices to show that for each $N \in M_K$, we have
\begin{equation}
    \label{suffices_thm1}
    \log P_N(\al) \ll - \psi(K).
\end{equation}
We apply Lemma \ref{upper_bound_key_lem} with ${\ell_0} = K$ and obtain
\[\log P_N(\al)  \leq \sum_{b = 1}^{b_{K-1} -1}\log \left\lvert 2\sin\left(\pi bq_{K-1}\delta_{K-1} 
    + \varepsilon_{K-1}(N)\right) \right\rvert
    + \mathcal{O}(K \log K).\]
Note that we have $\varepsilon_{K-1}(N) = 0$ and $b_{K-1}(N) \leq \left(\frac{1}{2}-\delta\right)a_K$, so
since $\log \lvert 2 \sin(\pi x)\rvert$ is monotonically increasing on $[0,1/2]$, we have for some $c = c(\delta) > 0$ that

\[ \sum_{b = 0}^{b_{K-1} -1}\log \left\lvert 2\sin\left(\pi bq_{K-1}\delta_{K-1} 
    + \varepsilon_{K-1}(N)\right) \right\rvert \leq 
    a_K \int_{1}^{b_{K-1}/a_K}\log \lvert 2 \sin(\pi x)\rvert \, \mathrm{d}x\leq 
    -c\cdot a_K  \ll -\psi(K),\]
    which completes the proof.

\subsection{Proof of Theorem \ref{thm2}}
By the proof of Theorem \ref{thm1}, we can deduce that
\begin{equation}\label{thm2_thm1}\limsup_{K \to \infty}\frac{\#\{0 \leq N \leq q_K: \log P_N(\al) \leq -2\psi(K)\}}{q_K} \geq \frac{1}{2}.\end{equation}
By the reflection principle \eqref{reflection_princ}, we see that at most one of the inequalities
\[\log P_N(\al) \leq -2\psi(K),\qquad \log P_{q_K-N-1}(\al) \leq -2\psi(K)\] can be fulfilled,
hence there is equality in \eqref{thm2_thm1}.
Applying the reflection principle a second time implies
\[\limsup_{K \to \infty}\frac{\#\{0 \leq N \leq q_K: \log P_N(\al) \geq \psi(K)\}}{q_K} \geq \frac{1}{2},\]
which finishes the proof of the first part of Theorem \ref{thm2}. \\

To show equality in the case where $\liminf_{k \to \infty} \psi(k)/(k \log k) \geq C$,
let $q_{K-1} \leq M < q_K$ be an arbitrary integer and let
$a_{{\ell_0}} = \max\limits_{\ell \leq K}a_{\ell}$.
We define the sets
\[M^+ := \left\{N \leq M: b_{{\ell_0}-1}(N) \geq \frac{a_{{\ell_0}}}{2}\right\},
\quad M^- := \left\{N \leq M: b_{{\ell_0}-1}(N) \leq \frac{a_{{\ell_0}}}{2}\right\}
\]
and the function
\begin{align}
f: M^+ &\to M^-\\
N = \sum_{\ell=0}^{K-1}b_{\ell}q_{\ell} &\mapsto \sum_{\ell=0}^{K-1}\tilde{b}_{\ell}q_{\ell}
\end{align}
with
$\sum\limits_{\ell=0}^{K-1}b_{\ell}q_{\ell}$ being the Ostrowski expansion of $N$ and

\begin{align*}\tilde{b}_{\ell} := \begin{cases} a_{{\ell_0}} - b_{{{\ell_0}}-1} &\text{ if } \ell = {\ell_0},\\ b_k &\text{ otherwise.}\end{cases}\end{align*}
It is straightforward to check that $f$ is well-defined and injective, hence $\lvert  M^- \rvert \geq \frac{M}{2}.$ For arbitrary $N \in M^-$, we apply Lemma \ref{upper_bound_key_lem}
to obtain 

\begin{equation}\label{thm2_lem_appl}
\log P_N(\al) \leq \sum_{b = 1}^{b_{{\ell_0}-1} -1}\log \left\lvert 2\sin\left(\pi bq_{{\ell_0}-1}\delta_{\ell_0 -1}
    + \varepsilon_{{\ell_0}-1}(N)\right) \right\rvert
    + \mathcal{O}(K \log K).\end{equation}
    By \eqref{in_part_eps}, it follows that
    \[0\leq bq_{{\ell_0}-1}\delta_{\ell_0 -1}
    + \varepsilon_{{\ell_0}-1}(N)\leq \frac{1}{2},\qquad b = 1, \ldots, b_{{\ell_0}-1}-1,\]
so each summand on the right-hand side of \eqref{thm2_lem_appl} is negative. Thus, for $N \geq \sqrt{M}$, we have
for almost every $\al$ that
\[\log P_N(\al) \ll K \log K
\ll \log N \log \log N.\]
Choosing $C$ sufficiently large, this shows $\log P_N(\al) \leq \psi(\log N)$ for $N \in M^- \cap \{\lceil \sqrt{M}\rceil,\ldots,M\}$,
and as

\[\limsup_{M \to \infty} \frac{\lvert M^- \cap \{1,\ldots, \lfloor \sqrt{M}\rfloor \}\rvert}{M} = 0,\]
the result follows.\\

\subsection*{Acknowledgements}
The author is grateful to Bence Borda for various comments on an earlier version of this paper.

\end{document}